\documentclass[11pt,leqno]{article}
\def\ds{\displaystyle}
\def\ac{\alpha c}
\def\m{{\bf m}}
\def\n{{\bf n}}
\def\x{{\bf x}}
\def\j{{\bf j}}
\def\k{{\bf k}}
\def\({\left(}
\def\){\right)}

\def\prodpk{\prod^p_{k=1}}
\def\prodpj{\prod^p_{j=1}}
\def\sumpj{\sum^p_{j=1}}
\def\sumj{\sum_{\bf j}}
\def\sumk{\sum_{\bf k}}
\def\g{\gamma}
\def\gg{\gamma\gamma}

\def\th{\theta}

\def\beq{\begin{equation}}
\def\eeq{\end{equation}}
\def\beqr{\begin{eqnarray}}
\def\eeqr{\end{eqnarray}}

\title{$q$-Analogues of Some Multivariable Biorthogonal Polynomials}

\author{George Gasper\footnote{Dept.\ of Mathematics, Northwestern 
University, Evanston, IL\ 60208.} ~and Mizan~Rahman\footnote{School 
of Mathematics and Statistics, Carleton University, Ottawa, ON, K1S 
5B6, Canada. This work was supported in part by the NSERC grant 
\#A6197.}}

\begin{document}
{}
\maketitle

\begin{abstract}
In 1989 M.V.~Tratnik found a pair of multivariable biorthogonal 
polynomials $P_{\n}(\x)$ and $\bar P_{\m}(\x),$ which is not necessarily
the complex conjugate of $P_{\m}(\x),$ such that
\[
\int^\infty_{-\infty} \cdots \int^\infty_{-\infty} w(\x)P_{\n}(\x) 
\bar P_{\m}(\x) \prod^p_{j=1} dx_j = \mu_{\n,\m} \delta_{N,M},
\]
where $\x=(x_1,\ldots,x_p)$, $\n=(n_1,\ldots,n_p)$, 
$\m=(m_1,\ldots,m_p)$, $N=\sum^p_{j=1} n_j$, $M=\sum^p_{j=1}m_j$, 
$\mu_{\n,\m}$ is the constant of biorthogonality (which Tratnik did 
not evaluate), 
\begin{eqnarray*}
w(\x)&=&
\Gamma(A-iX)\Gamma(B+iX)\left|
\frac{\Gamma(c+iX)\Gamma(d+iX)}{\Gamma(2iX)} \right|^2
\prod^p_{k=1} \Gamma(a_k+ix_k)\Gamma(b_k-ix_k),\\
&& X= \sum^p_{k=1} x_k, \quad A=\sum^p_{k=1} a_k, \quad B=\sum^p_{k=1}b_k,
\end{eqnarray*}
and the $a$'s, $b$'s, $x$'s, $c$ and $d$ are real.  In the $q$-case we find that the
appropriate weight function is a product of a
multivariable version of the integrand in the Askey-Roy integral
and of the Askey-Wilson weight function in a single variable that depends on 
$x_1,\ldots,x_p.$

In a related problem we find a discrete 2-variable Racah type 
biorthogonality:
\[
\sum^N_{x=0} \sum^N_{y=0} w_N(x,y)F_{m,n}(x,y)G_{m',n'}(x,y) = 
\nu_{m,n} \delta_{m,m'} \delta_{n,n'},
\]
where
\begin{eqnarray*}
w_N(x,y)&=& \frac{(\alpha q/\gamma\gamma', ~\gamma'/c,~ \alpha 
cq/\gamma';~q)_N}
                  {(\alpha q,~ 1/c,~ \alpha cq/\gamma\gamma';~q)_N} \\
&& \times
\frac
{
\(1-\frac{\gamma\gamma' q^{2x-N-1}} {\alpha c}\)(1-cq^{2y-N})\( 
\frac{\gamma\gamma' q^{-N-1}}{\alpha 
c},\gamma;q\)_x(cq^{-N},\gamma';q)_y}
{\(1-\frac{\gamma\gamma' q^{-N-1}}{\alpha 
c}\)(1-cq^{-N})\(q,\frac{\gamma'q^{-N}}{\alpha 
c};q\)_x\(q,\frac{cq^{1-N}}{\gamma'};q\)_y}
\\
&& \times \frac{
(1/c;q)_{x-y}(q^{-N};q)_{x+y}}{\(\frac{\gamma\gamma'}{\alpha c} ;q\)_{x-y}
\(\frac{\gamma\gamma' q^{-N}}{\alpha};q\)_{x+y}
              } \alpha^{-x}(\gamma')^{x-y},
\end{eqnarray*}
and $F_{m,n}(x,y)$, $G_{m',n'}(x,y)$ are certain bivariate extensions
of the $q$-Racah polynomials. 

\end{abstract}

\setcounter{equation}{0}
\section{Introduction}
Wilson polynomials [13], defined by
\beq
P_n(x)=(a+b)_n(a+c)_n(a+d)_n{\ }_4F_3\left[
\begin{array}{c}
-n,n+a+b+c+d-1,a-ix,a+ix\\
a+b,a+c,a+d
\end{array};1\right]
\label{eq1.1}
\eeq
satisfy an orthogonality relation on the real line
\beq
\int^\infty_{-\infty}
P_n(x)P_m(x)w(x)dx = h_n\delta_{n,m},
\label{eq1.2}
\eeq
where
\beq
w(x)=\left|
\frac{\Gamma(a+ix)\Gamma(b+ix)\Gamma(c+ix)\Gamma(d+ix)}{\Gamma(2ix)}
\right|^2
\label{eq1.3}
\eeq
is the positive weight function (under the assumption that $a,b,c,d$ 
are real or occur in complex conjugate pairs), and
\beq
~~~~~~~~h_n=4\pi n!\, (n\!+\!a\!+\!b\!+\!c\!+\!d-1)_n
\frac{\Gamma(n\!+\!a\!+\!b)\Gamma(n\!+\!a\!+\!c)\Gamma(n\!+\!a\!+\!d)
\Gamma(n\!+\!b\!+\!c)\Gamma(n\!+\!b\!+\!d)\Gamma(n\!+\!c\!+\!d)}
      {\Gamma(2n\!+\!a\!+\!b\!+\!c\!+\!d)}
\label{eq1.4}
\eeq
is the normalization constant.  By Whipple's transformation it is 
easy to see that $P_n(x)$ is symmetric in $a,b,c,d$, and that
\begin{eqnarray}
P_n(x)&=&(a+b)_n (c-ix)_n (d-ix)_n{\ }_4F_3\left[
\begin{array}{c}
-n,1-c-d-n,a+ix,b+ix\\
a+b,1-c-n+ix,1-d-n+ix
\end{array};1\right]\\
&=&(b+a)_n (c+ix)_n (d+ix)_n{\ }_4F_3\left[
\begin{array}{c}
-n,1-c-d-n,a-ix,b-ix\\
a+b,1-c-n-ix,1-d-n-ix
\end{array};1\right].\nonumber
\label{eq1.5}
\end{eqnarray}
Corresponding to each of these forms M.V.~Tratnik [10] introduced a 
multivariable polynomial:
\begin{eqnarray}
P_n(\x)&=&(A+c)_N (A+d)_N \prod^p_{k=1} (a_k+b_k)_{n_k}\\
&&\times \sum_{\j}\frac{(N+A+B+c+d-1)_J(A-iX)_J}{(A+c)_J(A+d)_J}
\prod^p_{k=1}
\frac{(-n_k)_{j_k}(a_k+ix_k)_{j_k}}{(a_k+b_k)_{j_k} j_k!},\nonumber
\label{eq1.6}
\end{eqnarray}
\begin{eqnarray}
\bar P_n(\x)&=&(B+c)_N (B+d)_N \prod^p_{k=1} (b_k+a_k)_{n_k}\\
&&\times \sum_{\j}\frac{(N+A+B+c+d-1)_J(B+iX)_J}{(B+c)_J(B+d)_J}
\prod^p_{k=1}
\frac{(-n_k)_{j_k}(b_k-ix_k)_{j_k}}{(b_k+a_k)_{j_k} j_k!},\nonumber
\label{eq1.7}
\end{eqnarray}
\begin{eqnarray}
Q_{\n}(\x)&=&(c-iX)_N (d-iX)_N \prod^p_{k=1} (a_k+b_k)_{n_k}\\
&&\times \sum_{\j}\frac{(1-c-d-N)_J(B+iX)_J}{(1-c-N+iX)_J(1-d-N+iX)_J}
\prod^p_{k=1}
\frac{(-n_k)_{j_k}(a_k+ix_k)_{j_k}}{(a_k+b_k)_{j_k} j_k!},\nonumber
\label{eq1.8}
\end{eqnarray}
\begin{eqnarray}
\bar Q_{\n}(\x)&=&(c+iX)_N (d+iX)_N \prod^p_{k=1} (b_k+a_k)_{n_k}\\
&&\times \sum_{\j}
\frac{(1-c-d-N)_J(A-iX)_J}{(1-c-N-iX)_J(1-d-N-iX)_J}
\prod^p_{k=1}
\frac{(-n_k)_{j_k}(b_k-ix_k)_{j_k}}{(b_k+a_k)_{j_k} j_k!},\nonumber
\label{eq1.9}
\end{eqnarray}
where $\x = (x_1,x_2,\ldots,x_p)$, $\n = (n_1,n_2,\ldots,n_p)$, $\j = 
(j_1,j_2,\ldots,j_p)$, and
$X=\sum^p_{k=1} x_k$, $N=\sum^p_{k=1} n_k$, $M=\sum^p_{k=1} m_k,$ $A=\sum^p_{k=1} a_k$, 
$B=\sum^p_{k=1} b_k$, $J=\sum^p_{k=1}j_k$, and the sums in 
(1.6)--(1.9) are from $j_k=0$ to $n_k$, $k=1,\ldots,p$.  Each of the polynomials in
(1.6)--(1.9) is
of (total) degree $2N$ in the variables $x_1,x_2,\ldots,x_p.$ The
overbars in (1.7), (1.9), and in (1.21) below are used to denote  
distinct systems of polynomials
and should not be confused with complex conjugation.
 Tratnik proved that
\beq
\int^\infty_{-\infty} \cdots \int^\infty_{-\infty} P_{\n}(\x)\bar 
P_{\m}(\x) w(\x) \prod^p_{k=1} dx_k=0, \hbox{ ~~if } N\not = M, 
\label{eq1.10}
\eeq
\beq
\int^\infty_{-\infty} \cdots \int^\infty_{-\infty}Q_{\n}(\x)\bar 
Q_{\m}(\x) w(\x) \prod^p_{k=1} dx_k=0, \hbox{ ~~if } N\not = M, 
\label{eq1.11}
\eeq
\beq
\int^\infty_{-\infty} \cdots \int^\infty_{-\infty}  P_{\n}(\x) 
Q_{\m}(\x) w(\x) \prod^p_{k=1} dx_k=0, \hbox{ ~~if } \n\not = \m,
\label{eq1.12}
\eeq
and
\beq
\int^\infty_{-\infty} \cdots \int^\infty_{-\infty} \bar 
P_{\n}(\x)\bar Q_{\m}(\x) w(\x) \prod^p_{k=1} dx_k=0, \hbox{ ~~if } 
\n\not = \m,
\label{eq1.13}
\eeq
where
\beq
w(\x)=\left| \frac{\Gamma(c+iX)\Gamma(d+iX)}{\Gamma(2iX)}\right|^2
\Gamma(A-iX)\Gamma(B+iX)\prod^p_{k=1}\Gamma(a+ix_k)\Gamma(b-ix_k).
\label{eq1.14}
\eeq
\noindent
Note that in (1.12) and (1.13) the biorthogonality holds in all of the indices
$n_1,n_2,\ldots,n_p,$ while in (1.10) and (1.11) the biorthogonality is  for
polynomials of different degrees ($N\ne M$).

Since Whipple's ${}_4F_3$ transformation does not apply for $p\geq 2$ 
the $P$'s and $Q$'s are no longer equivalent and hence the 
orthogonality in a single variable becomes biorthogonality in many 
variables.

We were curious to see what their $q$-analogues would be.  At first 
sight it might appear that they could be found in a pretty 
straightforward manner.  We were in for a surprise.  The first hurdle 
is an appropriate analogue of the weight function in (1.14).  There 
are many possible candidates but the one that works for a $q$-analogue of (1.10) is:
\begin{eqnarray}
w^{(p)}(\x;q) &:=& \frac{1}{(2\pi)^p}
\frac{(e^{2i\Theta},e^{-2i\Theta};q)_\infty}
      {(Ae^{-i\Theta},Be^{i\Theta};q)_\infty h(\cos \Theta;c,d;q)
\(\frac{\beta b_1}{B} e^{i\Theta}, \frac{qB}{\beta 
b_1}e^{-i\Theta};q\)_\infty}\\
&&\times \prodpk
\frac{(\beta_k e^{i\theta_k},q\beta^{-1}_k 
e^{-i\theta_k};q)_\infty}{(a_k 
e^{i\theta_k},b_ke^{-i\theta_k};q)_\infty}, \quad p\geq 2,\nonumber
\label{eq1.15}
\end{eqnarray}
where $-\pi\leq \theta_k \leq \pi$, $\theta_k=x_k \log q$ so that $e^{i\theta_k}=q^{i x_k}$ 
for $k=1,\ldots,p$,
$\Theta = \sumpj 
\theta_j$, $A=\prodpj a_j$, $B=\prodpj b_j$, 
 $h(\cos \Theta;c,d;q)$ is defined as in [2, (6.1.2)],$\, \beta$ is an arbitrary complex parameter such that
$\beta\neq  q^{\pm n}$ for $n=0$, $1,\ldots$, and
\beq
\beta_{k+1}=\frac{\beta_k}{a_k b_{k+1}} , \quad k=1,2,\ldots,p-1,
\label{eq1.16}
\eeq
with $\beta_1=\beta$.  By making repeated use of the Askey-Roy 
integral [2, (4.11.1)] followed by the use of the Askey-Wilson 
integral, we shall prove in section 2 that
\begin{eqnarray}
W^{(p)}(q) &:=& \int^\pi_{-\pi} \cdots  \int^\pi_{-\pi} w^{(p)}(\x;q) \prodpk d\theta_k\\
&=& \frac{2(ABcd;q)_\infty \prod^p_{k=2} (b_k\beta_k, q/b_k\beta_k ;q)_\infty}
{(q;q)^p_\infty (Ac,Ad,Bc,Bd,cd;q)_\infty \prodpk (a_kb_k;q)_\infty},\nonumber
\label{eq1.17}
\end{eqnarray}
which is also valid for $p=1$.  It is understood that the 
$(p-2)$-fold product in the numerator is taken to be 1 when $p=1$.

Let
\begin{eqnarray}
&&A_j = \prod^p_{k=j} a_k, ~B_j=\prod^p_{k=j}b_k, ~J_j=\sum^p_{k=j} 
j_k, ~K_j=\sum^p_{r=j}k_r,\\
&&N_j= \sum^p_{k=j}n_k, ~M_j=\sum^p_{k=j} m_k, ~\Theta_j = 
\sum^p_{k=j} \theta_k,\nonumber
\label{eq1.18}
\end{eqnarray}
so that
\beq
A_1=A, ~B_1=B, ~J_1=J, ~K_1=K, ~N_1=N, ~M_1=M, ~\Theta_1=\Theta.
\label{eq1.19}
\eeq
Analogous to Tratnik's polynomials in (1.6) and (1.7) we introduce 
the functions
\begin{eqnarray}
P_{\n}(\x;q)&=& (Ac,Ad;q)_N \prodpk (a_kb_k;q)_{n_k}\\
&&\times \sumj
\frac{(ABcdq^{N-1},Ae^{-i\Theta};q)_J}{(Ac,Ad;q)_J}q^J \prodpk
\frac{(q^{-n_k},a_ke^{i\theta_k};q)_{j_k}}{(q,a_kb_k;q)_{j_k}}\nonumber\\
&&\times
\frac{e^{i\left(j_1\Theta_{2}+\cdots+ j_{p-1}\Theta_p\right)}}
      {B_2^{j_1} \cdots B_p^{j_{p-1}}}
q^{-(N_2j_1+N_3j_2+\cdots+N_pj_{p-1})},\nonumber
\label{eq1.20}
\end{eqnarray}
and
\begin{eqnarray}
\bar P_{\m}(\x;q)&=& (Bc,Bd;q)_M \prodpk (a_k b_k;q)_{m_k}\\
&&\times \sumk
\frac{(ABcdq^{M-1},Be^{i\Theta};q)_K}{(Bc,Bd;q)_K} q^K\prod^p_{r=1}
\frac{(q^{-m_r},b_r e^{-i\theta_r};q)_{k_r}}{(q,a_r b_r;q)_{k_r}}\nonumber\\
&&\times
\frac{e^{i\(k_2(\Theta_2-\Theta)+\cdots+ k_{p}(\Theta_{p}-\Theta)\)}}
      {a_1^{K_2} a_2^{K_3} \cdots a_{p-1}^{K_p}}
q^{-\sum^p_{r=2}k_r(M-M_r)},\nonumber
\label{eq1.21}
\end{eqnarray}
Both $ P_{\n}(\x;q)$ and  $\bar P_{\m}(\x;q)$ are Laurent polynomials in the variables 
$q^{i x_1}, \ldots, q^{i x_p}.$ Note that if we divide $ P_{\n}(\x;q)$ by $(1-q)^{3N}$ and replace
its parameters $a_1,\ldots,a_p, b_1,\ldots,b_p, c,d,$  respectively, by
 $q^{a_1},\ldots,q^{a_p}, q^{b_1},\ldots,q^{b_p}, q^c,q^d,$  and then let $q\to 1,$
we obtain $ P_{\n}(\x)$ as a limit case.  Similarly, we see that $\bar P_{\m}(\x)$ is
limit case of $\bar P_{\m}(\x;q).$
In section 3 we shall do the integration and in section 4 prove the 
following $q$-analogue of (1.10):
\beq
P_{\n} \cdot \bar P_{\m} := \int^\pi_{-\pi} \cdots \int^\pi_{-\pi} 
P_{\n}(\x;q) \bar P_{\m}(\x;q)\, w^{(p)}(\x;q) \prodpk d\theta_k =0, \hbox{ ~~if } N\neq M,
\label{eq1.22}
\eeq
where $w^{(p)}(\x;q)$ is given by (1.15), and
\begin{eqnarray}
P_{\n} \cdot \bar P_{\m} &=&
L_p \sum^{m_1}_{k_1=0}\cdots \sum^{m_{p-1}}_{k_{p-1}=0}
q^{\sum^{p-1}_{j=1} k_j \sum^{j-1}_{r=0} (n_r-m_r)} \\
&&\times\,
\frac{\(ABcdq^{N-1},\frac{ABcdq^N}{a_pb_p};q\)_{k_1+\cdots+k_{p-1}}}
      {\(ABcdq^{N+m_p},\frac{ABcdq^{N-n_p}}{a_pb_p} ;q\)_{k_1+\cdots+k_{p-1}}}
\prod^{p-1}_{r=1}
\frac{(q^{-m_r},a_r b_r q^{n_r};q)_{k_r}}
      {(q,a_r b_r;q)_{k_r} },\nonumber
\label{eq1.23}
\end{eqnarray}
when $N=M$, with $n_0=1$ and $m_0=0$, and $L_p$ is as defined in (3.7).

\smallskip
Discrete multivariable extensions of the Racah polynomials were 
considered in Tratnik [12] as well as in van Diejen and Stokman [1]
and in Gustafson [5].
For other related works see, for instance, [4, 6, 9, 11].
We have found $q$-extensions of Tratnik's systems of multivariable Racah and Wilson 
polynomials, complete with their orthogonality relations, see this 
Proceedings [3] for our multivariable extension of the
Askey-Wilson polynomials.   However, there seems to be at least one more 
extension that, to our knowledge, has not yet been investigated.  The 
seed of this extension lies in Rosengren's [8] multivariable 
extension of the $q$-Hahn polynomials as well as in Rahman's [7] 
2-variable discrete biorthogonal system.  In sections 5 and 6 we 
shall prove the following 2-variable extension of the $q$-Racah 
polynomial orthogonality [2, (7.2.18)]:
\beq
\sum^N_{x=0} \sum^N_{y=0} 
w_N(x,y)F_{m,n}(x,y)G_{m',n'}(x,y)=\nu_{m,n}\delta_{m,m'} 
\delta_{n,n'},
\label{eq1.24}
\eeq
where $0\le m, n, m', n' \le N,$
\begin{eqnarray}
\lefteqn{~~~~~F_{m,n}(x,y)} \\
&=&
\frac{\(\frac{\alpha 
q^{N+1-x-y}}{\gg'};q\)_{m+n}(q^{x-y}/c;q)_n(\alpha 
cq^{1+y-x}/\gg';q)_m}
     {(q^{-N};q)_{m+n}(\alpha cq/\gg';q)_n(1/c;q)_m} 
c^{n-m}q^{mx+ny}\nonumber \\
&&\times \sum^m_{i=0} \sum^n_{j=0}
\frac{(q^{-m},\gamma q^x,\gg' q^{x-N-1}/\alpha 
c;q)_i(q^{-n},\gamma'q^y,cq^{y-N};q)_j(\gg'q^{-M-n}/\alpha;q)_{i+j} }
  {(q,\gamma,\gg'q^{x-y-m}/\alpha 
c;q)_i(q,\gamma',cq^{1+y-x-n};q)_j(\gg'q^{x+y-N-m-n}/\alpha;q)_{i+j} 
} q^{i+j},\nonumber
\label{1.25}
\end{eqnarray}
\begin{eqnarray}
\lefteqn{~~G_{m,n}(x,y)} \\
&=&  \sum^m_{i=0} \sum^n_{j=0}
\frac{(q^{-m},q^{-x},\gg'q^{x-N-1}/\alpha c;q)_i
(q^{-n},q^{-y},cq^{y-N};q)_j
(\alpha q^{m+n};q)_{i+j}}
     {(q,\g, \g'q^{n}/c;q)_i
(q,\gamma',\alpha cq^{m+1}/\gamma';q)_j
(q^{-N};q)_{i+j}} q^{i+j},
\nonumber
\label{1.26}
\end{eqnarray}
and the weight function is
\begin{eqnarray}
~~w_N(x,y)
&=& \frac{ (\alpha q/\gg',\gamma'/c, \alpha cq/\gamma';q)_N }
        { (\alpha q, 1/c, \alpha cq/\gg';q)_N } \\
&& \times
\frac{ (1-\gg'q^{2x-N-1}/\alpha c)(1-cq^{2y-N})(\gg'q^{-N-1}/\alpha 
c,\g;q)_x(cq^{-N},\g';q)_y }
{(1-\gg'q^{-N-1}/\alpha c)(1-cq^{-N})(q,\g'q^{-N}/\alpha 
c;q)_x(q,cq^{1-N}/\g';q)_y }\nonumber\\
&&\times
\frac{ (1/c;q)_{x-y}(q^{-N};q)_{x+y} }
{ (\gg'/\alpha c;q)_{x-y}(\gg'q^{-N}/\alpha;q)_{x+y} }
\alpha^{-x}(\g')^{x-y}.\nonumber
\label{1.27}
\end{eqnarray}
\noindent
The normalization constant in (1.24) is given by
\beq
~~~~~\nu_{m,n}=\frac{1-\alpha}{1-\alpha q^{2m+2n}}
\frac{(q,\alpha cq/\g';q)_m (q,\g'/c;q)_n (\alpha q/\gg',\alpha 
q^{N+1};q)_{m+n} }
      {(\g,1/c;q)_m (\g',\alpha cq/\gg';q)_n (\alpha,q^{-N};q)_{m+n} }
c^{n-m}q^{mn}.
\label{eq1.28}
\eeq
Notice that both $F_{m,n}(x,y)$ and $G_{m,n}(x,y)$ are Laurent polynomials 
in the variables $q^x$ and $q^y,$ and  $G_{m,n}(x,y)$
is a polynomial of (total) degree $n+m$ in the variables  $q^{-x} + \gg'q^{x-N-1}/\alpha c$ and
$q^{-y}+cq^{y-N}.$

We wish to make the observation that the summation in (1.24) is over 
the square of length $N$, although the vanishing of the weight 
function above the main diagonal, because of the factor 
$(q^{-N};q)_{x+y}$ in the numerator, makes it effectively over the 
triangle $0\leq x+y \leq N$.  A very innocuous observation but it 
will help simplify the calculations somewhat as we shall see in 
section 6.

It seems reasonable to expect that there is a multivariable extension 
of (1.24), but we were unable to find it, mainly because an extension 
of the $q$-shifted factorials of the type $(a;q)_{x-y}$ doesn't 
appear too obvious to us.

\setcounter{equation}{0}
\section{Calculation of $W^{(p)}(q)$}
The key to the proof of (1.17) is to observe that by periodicity we 
can change $\theta_1,\th_2,\ldots,\th_p$ to, say, 
$\Theta,\th_2,\ldots,\th_p$ (so that $\th_1=\Theta-\Theta_2)$, with 
the limits of integration unchanged.  So the total weight transforms 
to
\begin{eqnarray}
~~~W^{(p)}(q) &=&
\frac{1}{(2\pi)^{p-1}} \int^\pi_{-\pi}
\frac{(e^{2i\Theta},e^{-2i\Theta};q)_\infty d\Theta}
{(Ae^{-i\Theta},Be^{i\Theta};q)_\infty \, 
h(\cos\Theta;c,d;q)(\frac{\beta}{B_2}e^{i\Theta},\frac{qB_2}{\beta} 
e^{-i\Theta};q)_\infty}
\\
&&\times \int^\pi_{-\pi}\cdots \int^\pi_{-\pi} 
I_2(\theta_3,\ldots,\th_p) \prod^p_{k=3}
\frac{(\beta_k e^{i\th_k},qe^{-i\th_k}/\beta_k;q)_\infty}
      {(a_ke^{i\th_k}, b_ke^{-i\th_k};q)_\infty}
d\theta_3 \cdots d\th_p,\nonumber
\label{eq2.1}
\end{eqnarray}
where
\beq
I_2(\theta_3,\ldots,\th_p)=
\frac{1}{2\pi} \int^\pi_{-\pi}
\frac{(\beta_2 
e^{i\th_2},qe^{i(\Theta_3-\Theta)+i\th_2}/\beta,qe^{-i\th_2}/\beta_2,\beta 
e^{i(\Theta-\Theta_3)-i\th_2};q)_\infty}
{(a_2e^{i\th_2},b_1e^{i(\Theta_3-\Theta)+i\th_2},b_2e^{-i\th_2},a_1e^{i(\Theta-\Theta_3)-i\th_2};q)_\infty} 
d\th_2.
\label{eq2.2}
\eeq
However, this integral matches exactly with the Askey-Roy integral 
[2, (4.11.1)], provided we assume that $\max(|a_1|,|b_1|,|a_2|,|b_2|)<1$  (with, of course, $|q|<1$).  By [2,
(4.11.1)], it then  follows that
\beq
I_2(\theta_3,\ldots,\th_p)=
\frac{
(b_2\beta_2,q/b_2\beta_2, a_1a_2b_1b_2, a_1\beta_2 e^{i(\Theta - 
\Theta_3)}, qe^{i(\Theta_3-\Theta)}/a_1\beta_2; q)_\infty }
      {(q,a_1b_1,a_2b_2, a_1a_2 e^{i(\Theta-\Theta_3)}, 
b_1b_2e^{i(\Theta_3-\Theta)};q)_\infty} .
\label{eq2.3}
\eeq
Substitution of (2.3) into (2.1) makes it clear that the integration 
over $\theta_4$ presents exactly the same situation, and so does the 
remaining integrations up to and including $\th_p$.  Finally, one is 
left with an Askey-Wilson integral over $\Theta$:
\begin{eqnarray}
~~~W^{(p)}(q) &=&
\frac{(AB;q)_\infty \prod^p_{k=2}(b_k\beta_k,q/b_k\beta_k;q)_\infty}
      {(q;q)^{p-1}_\infty \prod^p_{k=1}(a_kb_k;q)_\infty}
\, \frac{1}{2\pi} \int^\pi_{-\pi}
\frac{(e^{2i\Theta},e^{-2i\Theta};q)_\infty }
{h(\cos\Theta;\, A,B,c,d;q)}\, d\Theta\\
&=&
\frac{2(ABcd;q)_\infty \prod^p_{k=2} (b_k\beta_k,q/b_k \, \beta_k;q)_\infty}
{(q;q)^p_\infty (Ac,Ad,Bc,Bd,cd;q)_\infty \prodpk (a_kb_k;q)_\infty} 
, \nonumber
\label{eq2.4}
\end{eqnarray}
by [2,(6.1.1)], which completes the proof of (1.17).

\setcounter{equation}{0}
\section{Computation of the integral in (1.22)}
We shall carry out the integrations in (1.22) in much the same way as 
we did in the previous section.  We transform the integration 
variables $\th_1,\ldots,\th_p$ to $\th_2,\ldots,\th_p$ and $\Theta$ 
as before; then we isolate the\break $\th_2$-integral by observing 
that the factors $~(a_1e^{i(\Theta-\Theta_3)-i\th_2} 
;q)_{j_1}\,(a_2e^{i\th_2};q)_{j_2}\,(b_1 e^{i(\Theta_3-\Theta)+i\th_2};q)_{k_1}$\break
$\times 
(b_2e^{-i\th_2};q)_{k_2}\,e^{i\th_2(j_1+k_2)+ik_2(\Theta_3-\Theta)+ij_1\Theta_3}$ 
can be glued on to the integrand of $W^{(p)}(q)$, to get
\[
(-\beta)^{j_1+k_2}q^{\({j_1+k_2 \atop 2}\)} e^{ij_1\Theta}
   \frac{1}{2\pi}\! \int^\pi_{-\pi}
\frac{\big(\beta_2e^{i\th_2},\, q^{1-j_1-k_2}\, 
e^{i(\Theta_3-\Theta)+i\th_2}/\beta,  \beta 
q^{j_1+k_2}e^{i(\Theta-\Theta_3)-i\th_2},
qe^{-i\theta_2}/\beta_2;q\big)_\infty}
{\big(a_2q^{j_2}e^{i\th_2},b_1q^{k_1}e^{i(\Theta_3-\Theta)+i\th_2},b_2 
q^{k_2}e^{-i\th_2},a_1q^{j_1}e^{i(\Theta-\Theta_3)-i\th_2};q\big)_\infty}d\th_2
\]
which via [2, (4.11.1)] equals, on a bit of simplification,
\begin{eqnarray}
&& a_1^{k_2} b_2^{j_1}q^{j_1 k_2}e^{ij_1\Theta_3}
\frac{\big(b_2\beta_2,q/b_2\beta_2,a_1a_2b_1b_2q^{j_1+j_2+k_1+k_2};q\big)_\infty}
{\big(q, a_1b_1q^{j_1+k_1},a_2b_2q^{j_2+k_2};q\big)_\infty}\\
&& \quad \times~ 
\frac{\big(a_1\beta_2e^{i(\Theta-\Theta_3)}, 
{qe^{i(\Theta_3-\Theta)}/a_1\beta_2};q\big)_\infty}
{\big(a_1a_2q^{j_1+j_2}e^{i(\Theta-\Theta_3)},b_1b_2q^{k_1+k_2}e^{i(\Theta_3-\Theta)};q\big)_\infty}.\nonumber
\label{eq3.1}
\end{eqnarray}
Since $\Theta_3 = \theta_3+\Theta_4$, we may now isolate the 
$\theta_3$--integral in exactly the same way, carry out a similar 
integration, simplify, and obtain
\begin{eqnarray}
\lefteqn{a_1^{k_2}(a_1a_2)^{k_3} (b_2b_3)^{j_1} b_3^{j_2} 
e^{i(j_1+j_2)\Theta_4}q^{j_1 k_2+(j_1+j_2)k_3}}\\
&\times~ 
\frac{\big(b_2\beta_2,\, q/b_2\beta_2,\, 
b_3\beta_3,q/b_3\beta_3,\,a_1a_2a_3b_1b_2b_3q^{j_1+j_2+j_3+k_1+k_2+k_3};q\big)_\infty}
{\big(q,q, 
a_1b_1q^{j_1+k_1},a_2b_2q^{j_2+k_2},a_3b_3q^{j_3+k_3};q\big)_\infty}\nonumber\\
&\times~ 
\frac{\big( a_1a_2\beta_3 e^{i(\Theta-\Theta_4)},\, 
{qe^{i(\Theta_4-\Theta)}/a_1a_2\beta_3};q\big)_\infty}
{\big(a_1a_2a_3q^{j_1+j_2+j_3}\, e^{i(\Theta-\Theta_4)},\, 
b_1b_2b_3q^{k_1+k_2+k_3}\, e^{i(\Theta_4-\Theta)};q\big)_\infty}.
\nonumber
\label{eq3.2}
\end{eqnarray}
A clear pattern is now emerging.  The $\th_p$ integral is
\begin{eqnarray}
&& ~~~q^{j_1k_2+(j_1+j_2)k_3+\cdots+(j_1+\cdots+j_{p-2})k_{p-1}}
 (a_1^{k_2+\cdots+k_{p-1}}
 a_2^{k_3+\cdots+k_{p-1}}
\cdots a_{p-2}^{k_{p-1}}) (b_2^{j_1}\,b_3^{j_1+j_2}\cdots 
b_{p-1}^{j_1+\cdots+j_{p-2}})\\
& \times &\! \!\! \!
\frac{\big(\frac{AB}{a_pb_p} q^{J+K-j_p-k_p};q\big)_\infty 
\ds\prod^{p-1}_{r=2} (b_r\beta_r,q/b_r\beta_r;q)_\infty}
{(q;q)^{p-2}_\infty \ds\prod^{p-1}_{r=1} (a_{r} b_{r} q^{j_{r} + 
k_{r}};q)_\infty }\nonumber\\
&\times &\! \!\! \!
\left[
\frac{e^{-i\Theta k_p}}{2\pi}
  \int^\pi_{-\pi}                  
\frac{
\left(\beta_p e^{i\th_p},
\frac{
qe^{i(\th_p-\Theta)}}{(a_1\cdots a_{p-2})\beta_{p-1}},
\, qe^{-i\th_p}/\beta_p,\, (a_1\cdots a_{p-2})\beta_{p-1} 
e^{i(\Theta-\th_p)};q\)_\infty
}
{
\(a_pq^{j_p}e^{i\th_p},\frac{B}{b_p} q^{K-k_p}e^{i(\th_p-\Theta)},b_p 
q^{k_p}e^{-i\th_p},
\frac{A}{a_p}
q^{J-j_p} e^{i(\Theta-\th_p)}; q
\right)_\infty
}
 e^{i\th_p(J-j_p+k_p)} d\th_p\right].\nonumber
\label{eq3.3}
\end{eqnarray}
The expression in [~~~] above can, once again, be computed by use of 
[2, (4.11.1)], and simplified to
\beq
\left(\frac{Aq^{J-j_p}}{a_p}\right)^{k_p}
\frac{b_p^{J-j_p}\(b_p\beta_p,q/b_p\beta_p,\frac{A\beta_pe^{i\Theta}}{a_p},
\frac{qa_p}{A\beta_p} e^{-i\Theta},ABq^{J+K};q\)_\infty
}
{\(q,a_pb_pq^{j_p+k_p},Aq^Je^{i\Theta},Bq^K 
e^{-i\Theta},\frac{AB}{a_pb_p}q^{J+K-j_p-k_p};q\)_\infty
}.
\label{eq3.4}
\eeq
Since, by repeated application of (1.16) we get $A\beta_p/a_p=\beta 
b_1/B$, the $\Theta$-integral simply becomes the Askey-Wilson integral
\begin{eqnarray}
&&\frac{1}{2\pi} \int^\pi_{-\pi}
\frac{\(e^{2i\Theta}, e^{-2i\Theta};q\)_\infty }
      {h(\cos \Theta; Aq^J, Bq^K, c, d;q)}\, d\Theta\\
&& \quad =
\frac{2(ABcd\, q^{J+K};q)_\infty}
      {(q,cd,ABq^{J+K},Acq^J,Adq^J,Bcq^K,Bdq^K;q)_\infty}.\nonumber
\label{eq3.5}
\end{eqnarray}
Collecting these results and substituting into the integral in 
(1.22), we find that
\begin{eqnarray}
P_{\n}\cdot \bar P_{\m} &=& L_p\sum_{\j} \sum_{\k}
\frac{(ABcdq^{N-1};q)_J(ABcdq^{M-1};q)_K}{(ABcd;q)_{J+K} } q^{J+K}\\
&&\times \prod^p_{r=1}
\frac{(q^{-n_r};q)_{j_r} (q^{-m_r};q)_{k_r}(a_rb_r;q)_{j_r+k_r}}{(q,a_r b_r ;q)_{j_r}(q,a_rb_r;q)_{k_r} }
\nonumber\\
&&\times~ q^{\sum^p_{s=2} [j_{s-1}(K_s-N_s)+k_s(M_s-M)]},\nonumber
\label{eq3.6}
\end{eqnarray}
where
\beq
L_p=(Ac,Ad;q)_N(Bc,Bd;q)_M W^{(p)}(q) \prod^p_{r=1} 
(a_rb_r;q)_{m_r}(a_rb_r;q)_{n_r}.
\label{eq3.7}
\eeq

\setcounter{equation}{0}
\section{Biorthogonality}\label{sect4}
The sum over $j_1$ and $k_1$ in (3.6) gives
\begin{eqnarray}
&&
\frac{(ABcdq^{N-1};q)_{J_2}(ABcdq^{M-1};q)_{K_2}}
      {(ABcd;q)_{J_2+K_2} }
q^{J_2+K_2}\\
&&\times\sum^{m_1}_{k_1=0}
\frac{(q^{-m_1},ABcdq^{M+K_2-1};q)_{k_1} }
      {(q,ABcdq^{J_2+K_2};q)_{k_1}} q^{k_1} {}_3\phi_2
\left[ \begin{array}{lll}
q^{-n_1},& ABcdq^{N+J_2-1},&a_1b_1q^{k_1}\\
&ABcdq^{J_2+K_2+k_1},&a_1b_1
\end{array}; q,q^{1+K_2-N_2} \right].\nonumber
\label{eq4.1}
\end{eqnarray}

\noindent
Since, by [2, (3.2.7)], the above ${}_3\phi_2$ equals
\[
~~~\frac{(ABcd;q)_{J_2+K_2+k_1}(q^{1+K_2-N};q)_{n_1} }
      {(ABcd;q)_{J_2+K_2+n_1}(q^{1+K_2-N})_{k_1} } \,
{}_3\phi_2\left[
\begin{array}{lll}
q^{-k_1}, & ABcdq^{N+J_2-1}, &a_1b_1q^{n_1}\\
&ABcdq^{n_1+J_2+K_2}, &a_1b_1
\end{array} ;q, q^{1+K-N} \right],
\]
we can now do the summation over $k_1$ via [2, (1.5.3)] to obtain that 
the expression in (4.1) reduces to
\begin{eqnarray}
&&
\frac{(ABcdq^{N-1};q)_{J_2} (ABcdq^{M-1};q)_{K_2} 
(ABcdq^{N+M_2-1};q)_{m_1} (q^{1+K_2-N};q)_{n_1} }
{(ABcd;q)_{n_1+J_2+K_2} (q^{1+K_2-N};q)_{m_1} }\\
&& \quad\times~
(-1)^{m_1} q^{({m_1\atop 2})+(1+K_2-N)m_1+J_2+K_2}\nonumber\\
&& \quad \times~
{}_4\phi_3\left[
\begin{array}{llll}
q^{-m_1},& a_1b_1q^{n_1}, & ABcdq^{N+J_2-1}, & ABcdq^{M+K_2-1}\\
& a_1b_1,& ABcdq^{N+M_2-1}, &ABcdq^{n_1+J_2+K_2}
\end{array} ; q,q \right].\nonumber
\label{eq4.2}
\end{eqnarray}
Note that the ${}_4\phi_3$ series is balanced.  Now, the sum over 
$j_2$ and $k_2$ gives
\begin{eqnarray}
&&\frac{(ABcdq^{N-1};q)_{J_3} (ABcdq^{M-1};q)_{K_3} 
(ABcdq^{N+M_2-1};q)_{m_1}(q^{1+K_3-N};q)_{n_1}}
{(ABcd;q)_{n_1+J_3+K_3} (q^{1+K_3-N};q)_{m_1}}\\
&& \quad \times~ (-1)^{m_1} q^{({m_1\atop 2})+(1+K_3-N)m_1+J_3+K_3}\nonumber\\
&& \quad\times~ \sum^{m_1}_{k_1=0}
\frac{(q^{-m_1},a_1b_1q^{n_1},ABcdq^{N+J_3-1},ABcdq^{M+K_3-1};q)_{k_1}}
{(q, a_1b_1, ABcdq^{n_1+J_3+K_3},ABcdq^{M_2+N-1};q)_{k_1}} q^{k_1}\nonumber\\
&& \quad \times~ \sum^{m_2}_{k_2=0}
\frac{(q^{-m_2},ABcdq^{M+K_3+k_1-1},q^{1+K_3-N_2};q)_{k_2}}
{(q,ABcdq^{n_1+J_3+K_3+k_1},q^{1+K_3-N+m_1};q)_{k_2}} q^{k_2}\nonumber\\
&& \quad\times~  {}_3\phi_2
\left[ \begin{array}{lll}
q^{-n_2}, & ABcdq^{N+J_3+k_1-1},& a_2b_2q^{k_2}\\
&ABcdq^{n_1+J_3+K}, &a_2b_2
\end{array} ; q,q^{1+K_3-N_3} \right].\nonumber
\label{eq4.3}
\end{eqnarray}

\noindent
As in the previous step we apply [2, (3.2.7)] to the ${}_3\phi_2$ 
series above, use [2, (1.5.3)] to do the $k_2$ sum and simplify the 
coefficients to reduce (4.3) to the following expression
\begin{eqnarray}
&&\frac{(ABcdq^{N-1};q)_{J_3} (ABcdq^{M-1};q)_{K_3} 
(ABcdq^{N+M_3-1};q)_{m_1+m_2}(q^{1+K_3-N};q)_{n_1+n_2}}
{(ABcd;q)_{n_1+n_2+J_3+K_3} (q^{1+K_3-N};q)_{m_1+m_2}}\\
&&\quad \times~(-1)^{m_1+m_2} q^{({m_1+m_2\atop 2}) +(1+K_3-N)(m_1+m_2)+J_3+K_3}
\nonumber\\
&&\quad \times~ \sum^{m_1}_{k_1=0} \sum^{m_2}_{k_2=0}
\frac{(q^{-m_1}, a_1b_1q^{n_1};q)_{k_1}(q^{-m_2},a_2b_2q^{n_2};q)_{k_2}}
{(q, a_1b_1;q)_{k_1}(q,a_2b_2;q)_{k_2}}
q^{(n_1-m_1)k_2}\nonumber\\
&&\quad \times~
\frac{(ABcdq^{N+J_3-1},ABcdq^{M+K_3-1};q)_{k_1+k_2}}
{(ABcdq^{n_1+n_2+J_3+K_3},ABcdq^{M_3+N-1};q)_{k_1+k_2}} q^{k_1+k_2}.\nonumber
\label{eq4.4}
\end{eqnarray}
A clear pattern of terms is now emerging, and by induction we find 
that at the ($p-1$)--th step the sum over 
$j_1,k_1,\ldots,j_{p-1},k_{p-1}$ in (3.6) equals
\begin{eqnarray}
&&\frac{(ABcdq^{N-1};q)_{J_p}(ABcdq^{M-1};q)_{K_p} 
(ABcdq^{N+M_p-1};q)_{M-m_p}(q^{1+K_p-N};q)_{N-n_p}}
{(ABcd;q)_{N-N_p+J_p+K_p}(q^{1+K_p-N};q)_{M-m_p}}\\
&& \times~ (-1)^{M-m_p} q^{({M-m_p \atop 
2})+(1+K_p-N)(M-m_p)+J_p+K_p}\nonumber\\
&& \times~ \sum_{k_1,\ldots,k_{p-1}}
\left[ \prod^{p-1}_{r=1}
\frac{(q^{-m_r},a_rb_rq^{n_r};q)_{k_r}}
{(q,a_rb_r;q)_{k_r}}
  \right]
\frac{(ABcdq^{N+J_p-1},ABcdq^{M+K_p-1};q)_{K-k_p}}
{(ABcdq^{N-n_p+J_p+K_p},ABcdq^{M_p+N-1};q)_{K-k_p}}\nonumber\\
&&\times~ q^{k_1+k_2(1+n_1-m_1)+\cdots + k_{p-1}(1+n_1+\cdots + n_{p-2}-m_1- 
\cdots - m_{p-2})}.\nonumber
\label{eq4.5}
\end{eqnarray}
Using (4.5) we obtain that the sum over $\j$ and $\k$ in (3.6) equals
\begin{eqnarray}
&&\frac{(ABcdq^{N+m_p-1};q)_{M-m_p}}
{(ABcd;q)_{N-n_p}}
(-1)^{M-m_p} q^{({M-m_p \atop 2})+(1-N)(M-m_p)}\\
&&\times \sum^{m_1}_{k_1=0} \cdots \sum^{m_{p-1}}_{k_{p-1}=0}
q^{k_1+k_2(1+n_1-m_1)+\cdots + k_{p-1}(1+n_1+\cdots + n_{p-2}-m_1- 
\cdots - m_{p-2})}\nonumber\\
&&\times \left[
\prod^{p-1}_{r=1}
\frac{(q^{-m_r},a_rb_rq^{n_r};q)_{k_r}}
{(q,a_rb_r;q)_{k_r}} \right]
\frac{(ABcdq^{N-1}, ABcdq^{M-1};q)_{k_1+\cdots +k_{p-1}} }
{(ABcdq^{N-n_p},ABcdq^{N+m_p-1};q)_{k_1+\cdots +k_{p-1}}} S_p,\nonumber
\label{eq4.6}
\end{eqnarray}
where
\begin{eqnarray}
S_p&=& \sum^{m_p}_{k_p=0}
\frac{(q^{-m_p},ABcdq^{M+k_1+\cdots 
+k_{p-1}-1};q)_{k_p}(q^{1+M-N-m_p+k_p};q)_\infty}
{(q, ABcdq^{N-n_p+k_1+\cdots 
+k_{p-1}};q)_{k_p}(q^{1+k_p-n_p};q)_\infty}q^{k_p} \\
&& \times~ {}_3\phi_2
\left[ \begin{array}{lll}
q^{-n_p}, & ABcdq^{N+k_1+\cdots + k_{p-1}-1}, & a_pb_pq^{k_p}\\
&ABcdq^{N-n_p+k_1+\cdots + k_{p-1}+k_p},& a_pb_p \end{array}
; ~q,q \right].\nonumber
\label{eq4.7}
\end{eqnarray}
Note that the ${}_3\phi_2$ series is balanced, so by [2, (II.12)] it has the sum
\[
\frac{(q^{1+k_p-n_p},
       \frac{ABcd}{a_pb_p}
      q^{N-n_p+k_1+\cdots + k_{p-1}};q)_{n_p}}
{(\frac{q^{1-n_p}}{a_pb_p}, ABcdq^{N-n_p+K};q)_{n_p}}.
\]
Hence,
\begin{eqnarray}
S_p&=&
\frac{\(
\frac{ABcd}{a_pb_p} q^{N-n_p+k_1+\cdots + k_{p-1}} ;q\)_{n_p} }
{\(\frac{q^{1-n_p}}{a_pb_p}, ABcdq^{N+k_1+\cdots +k_{p-1}-n_p};q\)_{n_p}}\\
&&\times~ \sum^{m_p}_{k_p=0}
\frac{\(q^{-m_p}, ABcdq^{M+k_1+\cdots + 
k_{p-1}-1};q\)_{k_p}\(q^{1+M-N-m_p+k_p};q\)_\infty}
{\( ABcdq^{N+k_1+\cdots + k_{p-1}};q\)_{k_p} (q;q)_\infty}
q^{k_p}.\nonumber
\label{eq4.8}
\end{eqnarray}
First, let us suppose that $N\geq M\geq 0$.  Then it is clear from 
the right side of (4.8) that $S_p$ is zero unless $k_p\geq N-M+m_p$, 
as well as $m_p\geq k_p$.  So, we must have
\beq
m_p+(N-M) \leq k_p \leq m_p.
\label{eq4.9}
\eeq
This is a contradiction unless $N=M$, and then $k_p=m_p$.  In that case
\beq
S_p=
q^{m_p}
\frac{\( \frac{ABcd}{a_pb_p}q^{N-n_p+k_1+\cdots + 
k_{p-1}};q\)_{n_p}\(q^{-m_p},ABcdq^{N+k_1+\cdots + 
k_{p-1}-1};q\)_{m_p}  }
{\( \frac{q^{1-n_p}}{a_pb_p};q\)_{n_p} \(ABcdq^{N+k_1+\cdots + 
k_{p-1}-n_p};q\)_{m_p+n_p}}.
\label{eq4.10}
\eeq
On the other hand, if $M\geq N \geq 0$ then
\beq
m_p-(M-N) \leq k_p \leq m_p.
\label{eq4.11}
\eeq
So we get
\begin{eqnarray}
S_p&=&
q^{m_p+N-M}
\frac{\( \frac{ABcd}{a_pb_p} q^{N-n_p+k_1+\cdots + k_{p-1}};q\)_{n_p}}
{\( \frac{q^{1-n_p}}{a_pb_p}, ABcdq^{N-n_p+k_1+\cdots + k_{p-1}} ;q\)_{n_p}}\\
&&\times~
\frac{\( q^{-m_p},ABcdq^{M+k_1+\cdots + k_{p-1}-1};q\)_{m_p+N-M}}
{\(ABcdq^{N+k_1+\cdots + k_{p-1}};q\)_{m_p+N-M}}\nonumber\\
&&\times~
{}_2\phi_1 \left[\begin{array}{ll}
q^{N-M}, & ABcdq^{N+m_p+k_1+\cdots k_{p-1}-1}\\
&ABcdq^{2N-M+m_p+k_1+\cdots + k_{p-1}}\end{array} ; q,q \right].\nonumber
\label{eq4.12}
\end{eqnarray}
However, the above ${}_2\phi_1$ equals
\beq
\frac{\( q^{1+N-M};q\)_{M-N}}
{(ABcdq^{2N-M+m_p+k_1+\cdots + k_{p-1}};q)_{M-N}}
(ABcdq^{N+m_p+k_1+\cdots + k_{p-1}-1})^{M-N},
\label{eq4.13}
\eeq
which vanishes unless $N=M$.  This completes the proof of (1.22).

\medskip
Also, with $N=M$, (3.6), (4.6) and (4.10) give
\begin{eqnarray}
P_{\n} \cdot \bar P_{\m} &=& L_p
\frac{(ABcdq^{N-1};q)_N \(\frac{a_pb_pq^{1-N}}{ABcd};q\)_{n_p}}
{(ABcd;q)_{N+m_p}(a_pb_p;q)_{n_p}}
(-1)^N q^{-({N\atop 2})-m_p-n_p}(ABcdq^N)^{n_p}\\
&&\times \sum^{m_1}_{k_1=0}\cdots\sum^{m_{p-1}}_{k_{p-1}=0}
q^{k_1+k_2(1+n_1-m_1)+\cdots+ k_{p-1}(1+n_1+\cdots + 
n_{p-2}-m_1-m_2-\cdots - m_{p-2})}
\nonumber\\
&&\times~
\frac{\(ABcdq^{N-1}, \frac{ABcdq^N}{a_pb_p} ;q\)_{k_1+\cdots + k_{p-1}} }
{\( ABcdq^{N+m_p}, \frac{ABcdq^{N-n_p}}{a_pb_p} ;q\)_{k_1+\cdots + k_{p-1}}}
\prod^{p-1}_{r=1}
\frac{\(q^{-m_r}, a_rb_rq^{n_r};q\)_{k_r}}
{(q,a_rb_r;q)_{k_r}},\nonumber
\label{eq4.14}
\end{eqnarray}
which is, of course, the same as (1.23).
By taking p=2, e.g., in which case the series on the right hand side of (4.14)
becomes a terminating balanced $ _4\phi_3$ series, it is easily seen that in general
the above inner product does not vanish when $N = M$ and ${\bf n} \ne {\bf m}.$

In closing this section we would like to point out that unlike the 
$q\to 1$ case that corresponds to the Tratnik biorthogonalities, the 
$q$-analogues of $P_{\n}\cdot Q_{\m}$, $P_{\n}\cdot \bar Q_{\m}$ or 
$Q_{\n} \cdot \bar Q_{\m}$ do not seem to work out the same way as 
$P_{\n}\cdot \bar P_{\m}$.

\setcounter{equation}{0}
\section{Transformations of $F_{m,n}(x,y)$ and $G_{m,n}(x,y)$}\label{sect5}
We shall now address the problem of proving the biorthogonality 
relation (1.24).  First of all, it is very simple to use [2, (II.20)] 
to prove that
\beq
\sum^N_{x=0} \sum^N_{y=0} W_N(x,y) = 1.
\label{eq5.1}
\eeq
The forms of $F_{m,n}(x,y)$ and $G_{m,n}(x,y)$ that turn out to be 
most convenient for the summations in (1.24) are as follows:
\begin{eqnarray}
F_{m,n}(x,y)&=&
\frac{\( \frac{\gg' q^{-N}}{\alpha};q\)_{x+y} \(\frac{\gg'}{\alpha 
c};q\)_{x-y}}
{(q^{-N};q)_{x+y} (c^{-1};q)_{x-y}}
\left( \frac{\alpha}{\gg'}\right)^x
\(\frac{\alpha q^{N+n+1}}{\gg'}\)^m q^{Nn}\\
&&\times~ \sum^x_{j=0} \sum^y_{k=0}
\frac{\(\frac{\gg'}{\alpha} q^{-m-n};q\)_{j+k} \(q^{-x}, 
\frac{\gg'}{\alpha c}q^{x-N-1},\g q^m; q\)_j  }
{
\( \frac{\gg'}{\alpha} q^{-N};q\)_{j+k} \(q,\g, \frac{\gg' 
q^{-n}}{\alpha c};q\)_j}\nonumber\\
&&\times~
\frac{\( q^{-y}, cq^{y-N},\g' q^n;q\)_k }
      {\(q,cq^{1-m},\g';q \)_k} q^{j+k},\nonumber
\label{eq5.2}
\end{eqnarray}
and
\begin{eqnarray}
G_{m,n}(x,y)&=&
\frac{
(\alpha q^{N+1};q)_{m+n} \(
\frac{\alpha cq}{\g'};q\)_m \(
\frac{\g'}{c};q\)_n  }
{(q^{-N};q)_{m+n} \( \frac{\g'}{c};q\)_m \(\frac{\alpha c q}{\g'};q\)_m}
\( \frac{\g' q^{-N-1}}{\alpha c}\)^m
\(\frac{cq^{-N}}{\g'}\)^n \\
&&\times \sum^m_{j=0} \sum^n_{k=0}
\frac{(\alpha q^{m+n};q)_{j+k}(q^{-m}, \g q^x, \frac{\alpha c}{\g'} 
q^{N-x+1};q)_j}
{(\alpha q^{N+1};q)_{j+k} (q,\g, \frac{ac}{\g'} q^{n+1};q)_j}\nonumber\\
&&\times ~
\frac{(q^{-n},\g' q^y, \frac{\g'q^{N-y}}{c};q)_k}
{(q, \g', \frac{\g' q^m}{c};q)_k}
q^{j+k}, \hbox{ ~~assuming~ } 0\leq m+n\leq N. \nonumber
\label{eq5.3}
\end{eqnarray}
Since
\begin{eqnarray*}
&&{}_4\phi_3\left[
\begin{array}{llll}
q^{-m}, & \alpha q^{j+m+n}, & q^{-x}, & \frac{\gg'}{\alpha c} q^{x-N-1}\\
& \g, & \frac{\g' q^n}{c}, & q^{j-N}
\end{array} ; q,q \right]\\
&=&
\frac{\( \frac{\alpha cq^{j+1}}{\g'},\alpha q^{N+n+1};q\)_m}
{\(\frac{\g'q^n}{c}, q^{j-N};q\)_m}
\( \frac{\g' q^{-N-1}}
         {\alpha c}\)^m
{}_4\phi_3\left[
\begin{array}{llll}
q^{-m}, & \alpha q^{j+m+n},& \g q^x, & \frac{\alpha c}{\g'} q^{N-x+1}\\
& \g,& \frac{\alpha cq^{j+1}}{\g'}, & \alpha q^{N+n+1}
\end{array}
;q,q \right]
\end{eqnarray*}
and
\begin{eqnarray*}
&&{}_4\phi_3\left[
\begin{array}{llll}
q^{-n}, & \alpha q^{i+m+n}, & q^{-y}, & c q^{y-N}\\
& \g', & \frac{\alpha c q^{i+1}}{\g'}, & q^{m-N}
\end{array} ; q,q \right]\\
&=&
\frac{\( \frac{\g'q^{m}}{c},\alpha q^{N+1+i};q\)_n}
{ \(\frac{\alpha cq^{i+1}}{\g'}, q^{m-N};q\)_n}
\( \frac{c q^{-N}}{\g'}\)^n
{}_4\phi_3\left[
\begin{array}{llll}
q^{-n}, & \alpha q^{m+n+i},& \g' q^y, & \frac{\g' q^{N-y}}{c} \\
& \g',& \frac{\g' q^{m}}{c}, & \alpha q^{N+1+i}
\end{array}
;q,q \right]
\end{eqnarray*}
by [2, (III.15)], (5.3) follows from (1.26) with a bit of simplification.

To derive (5.2) from (1.25) we need two applications of [2, (III.15)] 
on each of the two ${}_4\phi_3$ series involved in (1.25).  First
\begin{eqnarray}
&&
{}_4\phi_3\left[\begin{array}{llll}
q^{-m}, & \frac{\gg'}{\alpha c}q^{x-N-1}, & \g q^x, & \frac{\gg' 
q^{j-m-n}}{\alpha}\\
& \gamma, & \frac{\gg' q^{x-y-m}}{\alpha c}, & \frac{\gg'}{\alpha} 
q^{x+y-N-m-n+j}
\end{array} ; q,q \right]\\
&=&
\frac{(q^{y-N},c^{-1}q^{n-y-j};q)_m}
{\( \frac{\alpha c}{\gg'} q^{1+y-x}, \frac{\alpha q^{N-x-y+n+1-j}}{\gg'};q\)_m}
\(\frac{\alpha c q^{N-x+1}}{\gg'} \)^m\nonumber\\
&&\times~
{}_4\phi_3\left[
\begin{array}{llll}
q^{-x}, & \frac{\gg' q^{x-N-1}}{\alpha c}, & q^{-m}, & 
\frac{\alpha}{\g'} q^{m+n-j}\\
& \g, & q^{y-N}, & c^{-1} q^{n-y-j} \end{array}
;q,q \right]\nonumber\\
&=& \frac{\( \frac{\alpha c}{\gg'}q^{1+y-x+m},\frac{\alpha 
q^{N-x-y+1+m+n}}{\gg'};q\)_{x-m}}
{\(q^{y+m-N},c^{-1}q^{m+n-y};q\)_{x-m}}
\( \frac{\gg'}{\alpha c}q^{x-N-1}\)^{x-m}\nonumber\\
&&\times~
\frac{\(\frac{\gg'}{\alpha}q^{x+y-N-m-n},cq^{1+y-x-n};q\)_j}
{\(\frac{\gg'}{\alpha}q^{y-N-n},cq^{1+y-m-n};q\)_j}\nonumber\\
&&\times~
{}_4\phi_3\left[
\begin{array}{llll}
q^{-x},& \frac{\gg'}{\alpha c}q^{x-N-1}, &\g q^m, & 
\frac{\gg'q^{j-m-n}}{\alpha}\\
& \g, & \frac{\gg'q^{-y}}{\alpha c},& \frac{\gg'}{\alpha}q^{j+y-N-n}
\end{array};q,q \right].\nonumber
\label{eq5.4}
\end{eqnarray}
Substituted into (1.25) this leads to another balanced ${}_4\phi_3$ series:
\[
{}_4\phi_3\left[
\begin{array}{llll}
q^{-n}, & cq^{y-N},& \g' q^y, & \frac{\gg'}{\alpha} q^{i-m-n}\\
& \g',& cq^{y-m-n+1}, & \frac{\gg'}{\alpha} q^{i+y-N-n}
\end{array}
;q,q \right]
\]
which, when transformed twice in the same manner as in (5.4), leads to
\begin{eqnarray}
&&
\frac{\(c^{-1}q^{m+n-y}, \frac{\alpha q^{N+1-y+n}}{\gg'};q\)_{y-n}}
{\(q^{m+n-N},\frac{\alpha cq^{n+1}}{\gg'};q\)_{y-n}}
\(cq^{y-N}\)^{y-n} 
\frac{\( \frac{\gg'}{\alpha} q^{y-N-n}, \frac{\gg'}{\alpha c} q^{-y};q\)_i}
{\(\frac{\gg'}{\alpha} q^{-N},\frac{\gg'}{\alpha c} q^{-n};q\)_i}\\
&& ~~~\times~
{}_4\phi_3\left[
\begin{array}{llll}
q^{-y},& cq^{y-N}, & \g' q^n, & \frac{\gg'}{\alpha}q^{i-m-n}\\
& \g',\, cq^{1-m}, & \frac{\gg'}{\alpha} q^{i-N}
\end{array};q,q \right].\nonumber
\label{eq5.5}
\end{eqnarray}

After some simplifications (5.4) and (5.5) give (5.2).  Denoting the 
left hand side of (1.24) by \break $F_{m,n} \cdot G_{m',n'}$, it 
follows that
\begin{eqnarray}
F_{m,n}\cdot G_{m',n'} &=&
A_{m,n,m',n'} \sum^N_{x=0} \sum^N_{y=0}
\frac{\(1-\frac{\gg'}{\alpha c} q^{2x-N-1}\) (1-cq^{2y-N}) 
\(\frac{\gg'q^{-N-1}}{\alpha c}, \g;q\)_x}
{\( 1-\frac{\gg'}{\alpha c} q^{-N-1}\) (1-cq^{-N})\(q,\frac{\g' 
q^{-N}}{\alpha c};q\)_x}\\
&& \times~
\frac{(cq^{-N},\g';q)_y}
{\(q, \frac{cq^{1-N}}{\g'};q\)_y}
\g^{-x}(\g')^{-y}\nonumber\\
&&\times~
\sum_j \sum_k
\frac{\(\frac{\gg' 
q^{-m-n}}{\alpha};q\)_{j+k}\(q^{-x},\frac{\gg'}{\alpha c} 
q^{x-N-1},\g q^m;q\)_j}
{\(\frac{\gg' q^{-N}}{\alpha};q\)_{j+k}\(q, \g, 
\frac{\gg'q^{-n}}{\alpha c};q\)_j}\nonumber\\
&&\times~
\frac{(q^{-y},cq^{y-N},\g'q^n;q)_k}
{(q,cq^{1-m},\g';q)_k}
q^{j+k}\nonumber\\
&&\times \sum_r \sum_s
\frac{( \alpha q^{m'+n'};q)_{r+s}\(q^{-m'},\g q^x, \frac{\alpha 
cq^{N-x+1}}{\g'};q\)_r}
{\(\alpha q^{N+1};q)_{r+s}(q,\g,\frac{\alpha c}{\g'} q^{n'+1};q\)_r}\nonumber\\
&&\times~
\frac{\(q^{-n'},\g'q^y, \frac{\g'}{c}q^{N-y};q\)_s}
{\(q,\g',\frac{\g'q^{m'}}{c};q\)_s}
q^{r+s},\nonumber
\label{eq5.6}
\end{eqnarray}
where
\begin{eqnarray}
A_{m,n,m',n'}&=&
\frac{
(\alpha q/\gg', \g'/c, \alpha cq/\g';q)_N (\alpha 
q^{N+1};q)_{m'+n'}\(\frac{\g'}{c};q\)_{n'} \( \frac{\alpha 
cq}{\g'};q\)_{m'}}
{(\alpha q, 1/c, \alpha cq/\gg';q)_N (q^{-N};q)_{m'+n'} \(\frac{\alpha cq}{\g'};q\)_{n'}
\(\frac{\g'}{c};q\)_{m'}}
\\
&&\times~
\( \frac{\g' q^{-N-1}}{\alpha 
c}\)^{m'}\(cq^{-N}/{\g'}\)^{n'}\(\frac{\alpha q^{N+n+1}}{\gg'}\)^m 
q^{Nn}.\nonumber
\label{eq5.7}
\end{eqnarray}

\setcounter{equation}{0}
\section{Proof of (1.24)}\label{sect6}
Since each term in the weight function can be glued on nicely with 
the $x$ and $y$ dependent terms of the two double series in (5.6), 
the $x,y$-sum can be isolated as
\begin{eqnarray*}
\lefteqn{
{}_6W_5\(\frac{\gg'}{\ac}q^{2j-N-1};\g q^{r+j}, \frac{\gg'}{\ac} q^j, 
q^{j-N};q,\g^{-1}q^{-j-r}\)
}\\
&& \times~ {}_6W_5\( cq^{2k-N};\g' q^{s+k}, cq^{k+1}, q^{k-N};q,(\g' 
q^{k+s})^{-1}\) \nonumber\\
&=&
\frac{\( \frac{\gg'}{\ac} q^{2j-N}, \frac{q^{-N-r}}{\g};q\)_{N-j} 
\(cq^{2k-N+1},\frac{q^{-N-s}}{\g'};q\)_{N-k}}
{\(\frac{\g' q^{j-r-N}}{\ac}, q^{j-N};q\)_{N-j} 
\(\frac{cq^{1-N+k-s}}{\g'},q^{k-N};q\)_{N-k}},
\nonumber
\end{eqnarray*}
by [2, (II.21)].  The sum over $j,\, k,\, r,\, s$ in (5.6) now reduces to
\begin{eqnarray}
F_{m,n}\cdot  G_{m',n'} &=&
A_{m,n,m',n'} \frac{\(\g q, \g' q, \ac q/\gg',1/c;q\)_N}
                 {(q,q,\ac q/\g', \g'/c;q)_N}
\\
&& \times~
\sum_j \sum_k \sum_r \sum_s
\frac{ \( \frac{\gg' q^{-m-n}}{\alpha};q\)_{j+k}(\alpha q^{m'+n'};q)_{r+s}}
      { \( \frac{\gg' q^{-N}}  {\alpha};q\)_{j+k}(\alpha 
q^{N+1};q)_{r+s}} \nonumber\\
&&\times~
\frac{\(q^{-N}, \frac{\gg'}{\ac},\g q^m;q\)_j (q^{-N},\g'q^n,cq;q)_k 
\(q^{-m'},\frac{\ac q}{\g'},\g q^{N+1};q\)_r}
      {\(q,\g, \frac{\gg' q^{-n}}{\ac} ;q\)_j (q,\g',cq^{1-m};q)_k 
\(q,\g,\frac{\ac q^{n'+1}}{\g'};q\)_r}\nonumber\\
&&\times~
\frac{(q^{-n'},\g'/c,\g'q^{N+1};q)_s}
{\(q,\g',\frac{\g' q^{m'}}{c};q\)_s}
\frac{(\g;q)_{r+j} (\g';q)_{s+k}}{(\g q;q)_{r+j}(\g' q;q)_{s+k}} 
q^{j+k+r+s} . \nonumber
\label{eq6.1}
\end{eqnarray}
The sum over $j$ is a multiple of
\begin{eqnarray}
&&{}_5\phi_4\left[ \begin{array}{lllll}
q^{-N},& \g q^r, & \frac{\gg'}{\ac},& \g q^m, & \frac{\gg'}{\alpha} q^{k-m-n}
\\
& \g q^{r+1}, & \frac{\gg' q^{-n}}{\ac}, & \g, & \frac{\gg'}{\alpha} q^{k-N}
\end{array} ;q,q \right]\\
&& \quad =
\frac{(q;q)_N \(\frac{\g' q^{-n-r}}{\ac};q\)_n (q^{-r};q)_m 
\(\frac{\g' q^{k-N-r}}{\alpha} ;q\)_{N-m-n}}
{(\g q^{r+1};q)_N \(\frac{\gg' q^{-n}}{\ac};q\)_n (\g;q)_m 
\(\frac{\gg'}{\alpha} q^{k-N};q\)_{N-m-n}}
(\g q^r)^N, \nonumber
\label{eq6.2}
\end{eqnarray}
by [2, (1.9.10)].  Together with a similar expression for the sum 
over $k$ we now have
\begin{eqnarray}
\lefteqn{F_{m,n}\cdot G_{m',n'}} \\
&=&
A_{m,n,m',n'} \frac{\( \frac{\ac q}{\gg'},1/c;q\)_N}{(\ac q/\g', \g'/c;q)_N}
\sum_r \sum_s \frac{(\alpha q^{m'+n'};q)_{r+s}}{(\alpha 
q^{N+1};q)_{r+s}}\nonumber\\
&&\times~
\frac{\(q^{-m'},\frac{\ac q}{\g'},\g q^{N+1};q\)_r 
\(q^{-n'},\frac{\g'}{c}, \g'q^{N+1};q\)_s}
      {\( q,     \frac{\ac q^{n'+1}}{\g'};q\)_r \(q, \frac{\g' q^{m'}}{c};q\)_s}
q^{r+s} \nonumber\\
&&\times~
\frac{\(\frac{\g' q^{-n-r}}{\ac};q\)_n \(\frac{\g' 
q^{-N-r}}{\alpha};q\)_{N-m-n}(q^{-r};q)_m}
{(\g q^{N+1};q)_r \(\frac{\gg' q^{-n}}{\ac};q\)_n (\g;q)_m 
\(\frac{\gg'}{\alpha} q^{-N};q\)_{N-m-n}}
(\g q^r)^N \nonumber\\
&&\times~
\frac{\( \frac{cq^{1-m-s}}{\g'};q\)_m 
\(\frac{q^{-N-r-s}}{\alpha};q\)_{N-m-n}(q^{-s};q)_n}
{(\g'q^{N+1};q)_s (cq^{1-m};q)_m (\g';q)_n\(\frac{\g'}{\alpha} 
q^{-N-r};q\)_{N-m-n}}
(\g' q^s)^N \nonumber\\
&=& A_{m,n,m',n'} \frac{\(\frac{\ac q}{\gg'},1/c;q\)_N}{\(\frac{\ac 
q}{\g'},\frac{\g'}{c};q\)_N} (\gg')^N
\frac{\(\frac{\g' q^{-n}}{\ac};q\)_n 
\(\frac{cq^{1-m}}{\g'};q\)_m(q^{-N}/\alpha;q)_{N-m-n}}
{\( \g', \frac{\gg' q^{-n}}{\ac};q\)_n (\g, cq^{1-m};q)_m 
\(\frac{\gg'}{\alpha} q^{-N};q\)_{N-m-n}}
\nonumber\\
&&\times~
\sum_r \sum_s
\frac{(\alpha q^{m'+n'};q)_{r+s}}{(\alpha q^{m+n+1};q)_{r+s}}
\frac{\(q^{-m'},\frac{\ac q^{n+1}}{\g};q\)_r \(q^{-n'}, \frac{\g' 
q^m}{c};q\)_s(q^{-r};q)_m}
{\(q, \frac{\ac q^{n'+1}}{\g'};q\)_r \(q, \frac{\g'}{c} q^{m'};q\)_s}
\nonumber\\
&& \quad \times~
(q^{-s};q)_n q^{(m+1)r+(n+1)s}.\nonumber
\label{eq6.3}
\end{eqnarray}
The $r,s$ sum is
\begin{eqnarray}
&&(-1)^{m+n} q^{({m+1\atop 2})+({n+1\atop 2})}
\frac{(\alpha q^{m'+n'};q)_{m+n}\(q^{-m'},\frac{\ac 
q^{n+1}}{\g'};q\)_m \(q^{-n'},\frac{\g' q^m}{c};q\)_n}
{\(\alpha q^{m+n+1};q\)_{m+n}\(\frac{\ac q^{n'+1}}{\g'};q\)_m 
\(\frac{\g' q^{m'}}{c};q\)_n}
\\
&& \times~
\sum^{m'-m}_{r=0} \sum^{n'-n}_{s=0}
\frac{(\alpha q^{m+n+m'+n'};q)_{r+s} \(q^{m-m'},\frac{\ac 
q^{n+m+1}}{\g'};q\)_r \(q^{n-n'},\frac{\g' q^{n+m}}{c};q\)_s}
{(\alpha q^{2m+2n+1};q)_{r+s} \(q, \frac{\alpha 
cq^{m+n'+1}}{\g'};q\)_r \(q, \frac{\g' q^{n+m'}}{c};q\)_s}
q^{r+s}\nonumber
\label{eq6.4}
\end{eqnarray}
which vanishes unless $m' \geq m$ and $n' \geq n$.

The sum in (6.4), via [2, (II.12) and (II.6)], equals
\[
\frac{ \( q^{1+m-m'+n-n'}, \frac{\ac}{\g'} q^{m+n+1};q\)_{n'-n} 
(q^{1+m-m'};q)_{m'-m}}
      { \(\alpha q^{2m+2n+1}, 
\frac{cq^{1-m'-n'}}{\g'};q\)_{n'-n}(\alpha q^{2m+n+n'+1};q)_{m'-m}}
\(\alpha q^{m+n+m'+n'}\)^{m'-m}
\]
which vanishes unless $m' \leq m$ and $n'\leq n$.  Thus we must have
\beq
F_{m,n}\cdot G_{m',n'} = 0 \quad \hbox{ unless~~ } (m,n) = (m',n'), 
\hbox{ ~and then}
\label{eq6.5}
\eeq
\beq
F_{m,n}\cdot G_{m,n} = \frac{1-\alpha}{1-\alpha q^{2m+2n}}
\frac{\(q, \frac{\ac q}{\g'};q\)_m\(q, \frac{\g'}{c};q\)_n 
\(\frac{\alpha q}{\gg'},\alpha q^{N+1};q\)_{m+n}}
{(\g, 1/c;q)_m \(\g', \frac{\ac q}{\gg'};q\)_n (\alpha, q^{-N};q)_{m+n}}
c^{n-m}  q^{mn},\nonumber
\label{eq6.6}
\eeq
which completes the proof of (1.24) and (1.28).

  \medskip
It may be mentioned that there are other double series 
representations for $F_{m,n}(x,y)$ that one could use instead of 
(5.2) in the derivation of the biorthogonality relation (1.24) which 
do not contain the factor $1/(q^{-N};q)_{x+y}$ that cancels out the 
$(q^{-N};q)_{x+y}$ factor in the weight function, but the subsequent 
computations turn out to be quite tedious, while the final result is, 
of course, the same.\\



\begin{thebibliography}{99}
\bibitem{1.}
J.F.~van Diejen and J.V.~Stokman,
{\it Multivariable $q$-Racah polynomials},
Duke  Math. J.\ {\bf 91} (1998), 89--136.

\bibitem{2.}
 G.~Gasper and M.~Rahman,
{\it Basic Hypergeometric Series},
Cambridge University Press, 1990.

\bibitem{3.}
\underline{$\qquad$} and  \underline{$\qquad$}, {\it Some systems of 
multivariable orthogonal Askey-Wilson polynomials}, This Proceedings (2003).

\bibitem{4.}
Ya. I.~Granovski\u\i \, and A.S.~Zhedanov,
{\it `Twisted' Clebsch-Gordan coefficients for $SU_q(2)$},
J.~Phys.~A {\bf 25}, (1992), L1029--L1032.

\bibitem{5.}
R.A.~Gustafson,
{\it A Whipple's transformation for hypergeometric series in $U(n)$ 
and multivariable hypergeometric orthogonal polynomials},
SIAM J.~ Math.\ Anal.\ {\bf 18} (1987), 495--530.

\bibitem{6.}
H.T.~Koelink and J.~Van der Jeugt,
{\it Convolutions for orthogonal polynomials from Lie and quantum 
algebra representations},
SIAM J.\ Math.\ Anal.\ {\bf 29} (1998), 794--822.

\bibitem{7.}
M.~Rahman,
{\it Discrete orthogonal systems corresponding to Dirichlet distribution},
Utilitas Mathematica, {\bf 20}  (1981), 261--272.

\bibitem{8.}
H.~Rosengren,
{\it Multivariable $q$-Hahn polynomials as coupling coefficients for 
quantum algebra representations},
Int.\ J.\ Math.\ Sci.\ {\bf 28} (2001), 331--358.

\bibitem{9.}
M.V.~Tratnik,
{\it Multivariable biorthogonal Hahn polynomials},
J.~Math.\ Phys.\ {\bf 30} (1989), 627--634.

\bibitem{10.}
\underline{$\qquad$},
{\it Multivariable Wilson polynomials},
J.~Math.\ Phys.\ {\bf 30}  (1989), 2001--2011.

\bibitem{11.}
\underline{$\qquad$},
{\it Some multivariable orthogonal polynomials of the Askey 
tableau---continuous families},
J.\ Math.\ Phys.\ {\bf 32} (1991), 2065--2073.

\bibitem{12.}
\underline{$\qquad$},
{\it Some multivariable orthogonal polynomials of the Askey 
tableau---discrete families},
J.\ Math.\ Phys.\ {\bf 32} (1991),
2337--2342.

\bibitem{13.} J.A.~Wilson,
{\it Some hypergeometric orthogonal polynomials},
SIAM J.\ Math.\ Anal.\ {\bf 11} (1980), 690--701.
\end{thebibliography}
\end{document}